\documentstyle[11pt]{article}
\title{An example concerning Ohtsuki's invariant and the full 
$SO(3)$ quautum invariant
 \thanks{The first author is supported partially by the National Science
 Foundation
 of P. R. China, the second author is supported partially by NSF grant 
 DMS 9304580}}
\author{ Bang-He Li \& Tian-Jun Li}
\date{}
\begin{document}
\maketitle
\baselineskip 18pt
\begin{abstract}

Two lens spaces are given to show that Ohtsuki's $\tau$  for rational homology
spheres does not determine Kirby-Melvin's $\{\tau_r^{'}, r{\hskip .1cm}
\mbox{odd}\geq3\}$
 \end{abstract}
\vskip 36pt

By using partial Kirby-Melvin's quautum $SO(3)$ invariants $\{\tau_r^{'}(M),
r{\hskip .1cm}\mbox{odd prime} > |H_1(M,Z)| \}$ [KM], Ohtsuki [O]
defined a topological invariant
$$\tau(M)=\sum_{n=0}^{\infty}\lambda_n(M)(t-1)^n \in Q[[t-1]]$$
for rational homology 3-sphere $M$.
R.Lawrence [La] conjectured that
$$\lambda_n(M)\in \cases {Z, & if $|H_1(M,Z)|=1$\cr
Z[ {1\over 2}, {1\over {|H_1(M,Z)|}}], & if $|H_1(M,Z)| > 1$\cr}$$
and if $r$ is an odd prime which does not divide $|H_1(M,Z)|$, then
$\{|H_1(M,Z)|\}_r\tau^{'}_r(M)$ is the $r-$adic limit of the series
$$\sum_{n=0}^{\infty}\lambda_n^{\surd}(M)h^n$$
where $\{\cdot\}_r$ stands for the Jacobi symbol, and
$h=e^{{2\pi i}\over r}-1$.

Rozansky [R] has proved that this conjecture is true. So $\tau (M)$ and
$\{ \tau_r^{'},
r$ odd prime not dividing $ |H_1(M,Z)| \}$ determine each other.

A natural question arises: Does $\tau$ determine all $\{\tau_r^{'}, r{\hskip .1cm}
\mbox{odd}\geq3\}$?

It was proved in [Li] that $\tau_r^{'}(M)=\tau_r^{'}(M')$ iff
$\xi_r(M,A)=\xi_r(M',A)$  for $r$ odd
$\geq 3$, where $A$ is any $r-$th primitive root of unit.
So the question is equivalent to: Does $\tau (M)$ determine all
$\xi_r(M,e_r)$? Where, $e_a$ stands for $e^{{2\pi i}\over a}$.

For lens space $L(p,q)$, all $\xi_r(L(p,q),e_r)$ has been
obtained in [LL1] ( explicit formulas for $\tau^{'}_r(L(p,q))$ were given in
[LL2]), that is:
let $r\geq 3$ be odd and $c=(p,r)$ the common factor, then
$$(1)\;\;\,\xi_r(L(p,q), e_r)=\cases{\displaystyle\{p\}_r e_{r}^{-12s(q,p)}
e_p^{r'(q+q^*)}{e_r^{2p^{'}}-e_r^{-2p^{'}}\over e_r^{2}-e_r^{-2}}\;\,,
&if $c=1$ \cr
\quad\cr
\quad\cr
\displaystyle (-1)^{{r-1\over2}{c-1\over2}}\{{p/ c}\}_{r/ c}
\{q\}_c
e_{r}^{-12s(q,p)}
\cr
\displaystyle e_{pc}^{({r/c})^{'}(q+q^*-\eta p^*p)}
e_{rc}^{-2\eta({p/c})^{'}}
{{\epsilon(c)\sqrt c\eta}\over {e_r^{-2}-e_r^{2}}} \;\,,
&if  $c>1$ ,  $c\mid q^*+\eta$  \cr
%\quad& and $r\equiv \pm 1 \pmod 4$\cr
\quad\cr
0,&if $c>1$ and $c\mid\!\llap /  q^{*}\pm 1$\cr}$$
where $\eta=1$ or $-1$, $p^*p+q^*q=1$ with $0< q^{*} <p$,
 $({p/ c})^{'}{p/ c}
+ ({r}/{c})'{{r}/{c}}= 1$, $s(p,q)$ is the Dedekind sum, and
$$\epsilon (c)=\cases {1, & if $c\equiv 1 \pmod 4$\cr
i, & if $c\equiv -1 \pmod 4$\cr}$$
\vspace{.5cm}
Since
$$\tau(L(p,q))=t^{-3s(q,p)}{{t^{1\over 2p}-t^{-{1\over 2p}}}\over
{t^{1\over 2}-t^{-{1\over 2}}}}$$
([O] and [LL2]), $\tau(L(p_1,q_1))=\tau(L(p_2,q_2))$ iff $p_1=p_2$ and
$s(q_1,p_1)=s(q_2,p_2)$. The following Theorem answers the question above.

{\bf Theorem.} $s(6,25)=s(11,25)$, while $\xi_r( L(25,11))=\xi_r( L(25,6))$
if and only if $(r,25)\not=5$.

Proof. We calculate $s(q,p)$ by the formula in [H]:
$$12s(q,p)=\sum_{i=1}^n m_i +{{q+q^*}\over p} -3n$$
if ${\displaystyle p\over q}=m_n-{{1}\over\displaystyle m_{n-1}-\cdots -{\strut{1}
\over\displaystyle m_2-{\strut{1}\over\displaystyle {m_1}}}}$ with
$m_i\geq 2$.

Now
$${25\over 6}=5-{1
\over\displaystyle 2-{\strut{1}
\over\displaystyle 2-{\strut{1}
\over\displaystyle 2-{\strut{1}
\over\displaystyle 2-\frac{1}{2}}}}},\;\;\;
{25\over 11}=3-{{1}
\over\displaystyle 2-{\strut{1}
\over\displaystyle 2-{\strut{1}
\over\displaystyle 3-\frac{1}{2}}}}$$
${11}^*=16, 6^*=21$, so
$$12s(6,25)=12 s(11,25)=-3+{27\over 25}$$
$c=(r,25)$ can be only $1, 5$ or $25$. If $c=25$, then $c\mid\!\llap / 6^*\pm 1$ and
$c\mid\!\llap / 11^*\pm 1$, thus by (1) ${\xi_r (L(25,11), e_r)}=
{\xi_r (L(25,6), e_r)}=0$. If $c=1$, it is easy to see from (1) that
the two $\xi_r$ are equal.

Assume $c=5$, then $c| 6^*-1$ and $11^*-1$. Now since $q^*q+p^*p=1$, we
have
$$q+q^*+(25)^*\times 25=\cases {27+125,& if $q=6$\cr
27+175, & if $q=11$\cr}$$
Therefore, by (1)
$$\frac{\xi_r (L(25,11), e_r)}{\xi_r (L(25,6), e_r)}=e_{125}
^{(r/ 5)^{'}\times 50}$$
Since  $({p/ c})^{'}{p/ c} + ({r}/{c})'{{r}/{c}}= 1$, we see that
$(r/ 5)^{'}$ is prime to $(p/ 5)=5$. This shows that $e_{125}
^{(r/ 5)^{'}\times 50}\not= 1$, and the theorem is proved.

 \begin{center} {\bf\Large \bf   References}\end{center}
 \begin{description}
\item{} [RT] Reshetikhin, N.Yu., Turaev, V.G.: { Invariants of 3-manifolds
via link polynomials and quantum groups}, Invent. Math. {\bf 103} (1991), 
547-597
\item{} [KM]    Kirby, R., Melvin, P.: { The 3-manifold invariants of Witten
and Reshetikhin-Turaev for $sl(2,C)$}, Invent. Math. {\bf 105} (1991), 473-545
\item{} [Li] Li, B.H., { Relations among Chern-Simons-Witten-Jones invariants
}, Science in China, series A, {\bf 38} (1995), 129-146
\item{} [LL1] Li, B.H., Li, T.J.:{ Generelized Gaussian Sums and Chern-Simons-
Witten-Jones invariants of Lens spaces}, J. Knot theory and its Ramifications,
vol.5 No. 2 (1996) 183-224
\item{} [O] Ohtsuki, T., { A polynomial invariant of rational homology
3-spheres}, Inv. Math. {\bf 123} (1996), 241-257
\item{} [H] Hickerson, D.,{ Continued fraction and density results}, J. Reine.
Angew. Math. {\bf 290} (1977), 113-116
\item{} [R] Rozensky, L.: { On p-adic properties of the Witten-
Reshetikhin-Turaev invariants}, math. QA/9806075, 12 Jun.1998
\item{} [La]   Lawrence, R.: { Asymptotic Expansions of Witten-
Reshetikhin-Turaev invariants for some simple 3-manifolds}, J. Math. Phys. {\bf
36} (1995) 6106-6129
\item{} [W] Witten, E: { Quatum field theory and the Jones polynomial},
Comm.Math. Phys., {\bf 121} (1989), 351-399
\end{description}
\vspace{.3cm}

 Author's address
$$\begin{tabular}{ll}
\mbox{ Bang-He Li                   }& \mbox{ Tian-Jun Li } \\
            \mbox{ Institute of Systems Science,  }&
    \mbox{ School of math. } \\
            \mbox{ Academia Sinica }&
    \mbox{ IAS } \\
    \mbox { Beijnig 100080}& \mbox{ Princeton NJ 08540 }  \\
    \mbox { P. R. China }& \mbox{ U.S.A.} \\
    \mbox{ Libh@iss06.iss.ac.cn} &
    \mbox{ Tjli@IAS.edu}\\
\end {tabular}$$
\end{document}